\documentclass[12pt,a4paper,draft]{article}
\usepackage[english]{babel}
\usepackage{latexsym,amsfonts,amssymb,amsmath,longtable,amsthm}
\renewcommand{\ge}{\geqslant}

\begin{document}
\newtheorem{lem}{{\bfseries Lemma}}
\newtheorem*{ttt}{{\bfseries Theorem}}
\newtheorem{cor}[lem]{{\bfseries Corollary}}
\newtheorem{prp}[lem]{{\bfseries Proposition}}
\theoremstyle{remark}
\newtheorem*{rem}{{\bfseries Remark}}
\newtheorem*{quest}{{\bfseries Problem}}
\theoremstyle{definition}
\newtheorem{tab}[lem]{{\bfseries Table}}
\newtheorem{df}[lem]{\bfseries Definition}
\newtheorem{con}[lem]{\slshape Conjecture}

\newcommand{\Sym}{\mathrm{Sym}}
\newcommand{\Alt}{\mathrm{Alt}}
\newcommand{\F}{\mathbb{F}}
\newcommand{\Z}{\mathbb{Z}}
\newcommand{\ov}{\overline}
\newcommand{\Aut}{\mathrm{Aut} }
\renewcommand{\P}{\mathbf{P}}
\newcommand{\la}{\langle}
\newcommand{\ra}{\rangle}

\begin{center}
{\Large\bfseries\slshape  ON THE EXISTENCE OF CARTER SUBGROUPS}

E.P.Vdovin\footnote{The work is supported by Russian Fond of Basic Research
(project 05--01--00797), grant of President of RF (МК--1455.2005.1) and SB RAS
(grant N~29 for young researches and Integration project~2006.1.2).}
\end{center}

{\small In the paper we obtain the existence criterion of a Carter subgroup in
a finite group in terms of its normal series. An example showing that the
criterion cannot be reformulated in terms of composition factors is given.}

\section{Introduction}

Recall that a nilpotent self-normalizing subgroup of a group $G$ is called a
{\em Carter subgroup}. The classical result by Carter \cite{Ca1} states that
every finite solvable group contains  Carter subgroups and all of
them are conjugate. A finite group $G$ is said to satisfy condition
{\bfseries(C)} if, for every its nonabelian composition factor $S$ and
for every its nilpotent subgroup $N$, Carter subgroups (if exist) of $\langle
\Aut_N(S),S\rangle$ are conjugate (the definition of $\Aut_N(S)$ is given
below). In the recent paper \cite[Theorem~10.1]{VdoAlmSimp} it is
proven that in every almost simple group with known simple socle Carter
subgroups are conjugate. Thus, modulo the classification of finite simple
groups, 
in every finite group Carter subgroups are conjugate. In the paper by a finite
group we almost mean a finite group satisfying {\bfseries(C)}, thus the results
of the paper does not depend on the classification of finite simple groups.
There exist finite
groups without Carter subgroups, the minimal  example
is~$\mathrm{Alt}_5$. In the paper we give a criterion of existence of Carter
subgroups in terms of normal series. 

If $G$ is a group, $A,B,H$ are subgroups of $G$ and $B$ is normal in $A$
($B\unlhd A$), then $N_H(A/B)=N_H(A)\cap N_H(B)$. If $x\in N_H(A/B)$, then $x$
induces an automorphism $Ba\mapsto B x^{-1}ax$ of $A/B$. Thus, there is a
homomorphism of $N_H(A/B)$ into $\mathrm{Aut}(A/B)$. The image of this
homomorphism is
denoted by $\mathrm{Aut}_H(A/B)$ while its kernel is denoted by $C_H(A/B)$. In
particular, if $S$ is a composition factor of $G$, then for any $H\leq G$ the
group $\mathrm{Aut}_H(S)$ is defined. 

Let
$G=G_0\geq G_1\geq\ldots\geq G_n=\{e\}$ be a chief series of $G$ (recall that
$G$ is assumed to satisfy {\bfseries(C)}). Then
$G_i/G_{i+1}=T_{i,1}\times\ldots\times T_{i,k_i}$, where
$T_{i,1}\simeq\ldots\simeq T_{i,k_i}\simeq T_i$ and $T_i$ is a simple group. 
If $i\ge 1$, then denote by $\ov{K}_i$ a Carter subgroup of $G/G_i$ (if it
exists) and by $K_i$
its complete preimage in $G/G_{i+1}$. If $i=0$, then $\ov{K}_0=\{e\}$ and
$K_0=G/G_1$.  We say
that a finite group $G$ satisfies condition {\bfseries (E)}, if, for every
$i,j$, either $\ov{K}_i$ does not exist, or 
$\Aut_{K_i}(T_{i,j})$ contains a Carter subgroup. 

The following lemma shows that the homomorphic image of a Carter subgroup is a
Carter subgroup. We shall use this fact substantially.

\begin{lem}\label{InhByHomomorphism}{\em \cite[Lemma~4]{Vdoconj}}
Let $G$ be a finite group, $H$ be a normal
subgroup of $G$ satisfying {\em\bfseries(C)}, and $K$ be a Carter subgroup of
$G$. Then $KH/H$ is
a Carter subgroup of~${G/H}$.
\end{lem}

\begin{proof}
By \cite{Vdoconj}, Carter subgroups of $KH$ are conjugate. Assume that there
exists $x\in N_{G}(KH)$. Then $K^x$ is a Carter subgroup of $KH$. Since Carter
subgroups of $KH$ are conjugate, there exists $y\in KH$ such that
$K^x=K^{y^{-1}}$, so $xy\in N_G(K)$. Since $K$ is a Carter subgroup of 
$G$, we obtain that
$xy\in K$ and $x\in K\leq KH$, a contradiction.
\end{proof}

\section{Criterion}

\begin{lem}\label{CarterSubgroupsInCompFactorsUnderHomomorphism}
Let $G$ be a finite group, $H$ be a normal subgroup of $G$ and $S$ be a
composition factor of~$G/H$ {\em (}hence of $G$ as well{\em )}.

Then $\Aut_G(S)=\Aut_{G/H}(S)$.
\end{lem}

\begin{proof}
Note that there exists a surjective homomorphism
$\varphi:\Aut_G(S)\rightarrow\Aut_{\ov{G}}(S)$, defined by
$$\Aut_G(S)=\left(N_G(A)\cap N_G(B)\right)/C_G(A/B)\rightarrow \left(N_G(A)\cap
N_G(B)\right)/\left(HC_G(A/B)\right)=\Aut_{\ov{G}}(S).$$
Since $F^{\ast}(\Aut_G(S))=S\cap \mathrm{Ker}(\varphi)=\{e\}$ it folows that
$\mathrm{Ker}(\varphi)=\{e\}$.
\end{proof}

Below we shall need to know some additional information about the structure of
Carter subgroups in groups of special type. Let $A'$ be a group with a normal
subgroup $T'$. Consider the direct product
$A_1\times\ldots\times A_k$, where $A_1\simeq\ldots\simeq A_k\simeq A'$ and its
normal subgroup $T=T_1\times\ldots\times T_k$, where $T_1\simeq \ldots\simeq
T_k\simeq T'$.
Consider the symmetric group
$\Sym_k$, acting on $A_1\times\ldots\times A_k$ by $A_i^s=A_{i^s}$ for all
$s\in S$ and define $X$ to be a semidirect product $\left(A_1\times\ldots\times
A_k\right)\leftthreetimes \Sym_k$ (permutation wreath product of
$A'$ and $\Sym_k$). Denote by $A$ the direct product
$A_1\times\ldots\times A_k$ and by $\pi_i$ the projection
$\pi_i:A\rightarrow A_i$.  In the introduced notations
the following lemma holds.

\begin{lem}\label{CarterInSubDirectProductOfSoluble}
Let $G$ be a subgroup of $X$ such that $T\leq G$,  $G/(G\cap T)$ is nilpotent
and $(G\cap A)^{\pi_i}=A_i$. Assume also that $A$ is solvable. Let $K$ be a
Carter subgroup of~$G$.

Then $(K\cap A)^{\pi_i}$ is a Carter subgroup of~$A_i$.
\end{lem}

\begin{proof}
Assume that the statement is not true and let $G$ be a counterexample of minimal
order with minimal $k$. Then $S=G/(G\cap A)$ is
transitive and
primitive. Indeed, if $S$ is not transitive, then $S\leq
\Sym_{k_1}\times\Sym_{k-k_1}$, hence $G\leq G_1\times G_2$. If we denote  by
$\psi_i:G\rightarrow G_i$ the natural homomorphism, then $G^{\psi_i}=G_i$
satisfies conditions of the lemma and $K^{\psi_i}=K_i$ is a Carter subgroup of
 $G_i$. Clearly $(G\cap A)^{\pi_j}=(G_i\cap A^{\psi_i})^{\pi_j}$,
where $i=1$ if $j\in\{1,\ldots,k_1\}$ and $i=2$ if $j\in\{k_1+1,\ldots,k\}$,
thus we obtain the statement by induction. If $S$ is transitive, but is
not primitive, let $\Omega_1=\{T_1,\ldots,T_m\},
\Omega_2=\{T_{m+1},\ldots,T_{2m}\}, \ldots,
\Omega_l=\{T_{(l-1)m+1},\ldots,T_{lm}\}$ be a system of imprimitivity. Then it
contains a nontrivial intransitive normal subgroup
$$F'\leq
\underbrace{\Sym_{m}\times\ldots\times \Sym_{m}}_{l\text{ times}},$$ where
$k=m\cdot l$. Consider the complete preimage $F$ of $F'$ in $X$. Then
$G\cap F\leq F_1\times\ldots\times F_l$. Denote by
$\psi_i:F\rightarrow F_i$ the natural projection, then $(G\cap
F)^{\psi_i}=F_i$. Note that all of $F_i$ satisfy conditions of the lemma and, if
we define $T_i'=T_{(i-1)m+1}\times\ldots\times T_{im}$, then $G$ satisfies
conditions of the lemma with $T'=T'_1\times\ldots\times T'_l$ and $A'=F$. By
induction we have that $(K\cap F)^{\psi_i}$ is a
Carter subgroup of $F_i$ and, if $j\in\{m\cdot(i-1)+1,\ldots m\cdot i\}$,
then $\left((K\cap F)^{\psi_i}\cap A^{\psi_i}\right)^{\pi_j}$ is a Carter
subgroup of $A_j$. Since $(G\cap A)^{\pi_j}=\left((K\cap
F)^{\psi_i}\cap A^{\psi_i}\right)^{\pi_j}$ (for suitable $i$), we obtain
the statement by induction.

Let $Y'$ be a minimal normal subgroup of $G$ contained in
$T$ (if $Y'$ is trivial, then $T$ is trivial and we have nothing to prove,
since $G$ is nilpotent in this case). Thus $Y'$ is a normal elementary Abelian
$p$-group. Let
$Y_i=(Y')^{\pi_i}$, then $Y=Y_1\times\ldots\times Y_k$ is a nontrivial normal
subgroup of
$G$ ($Y$ is a subgroup of $G$ since $T\leq G$). 
Let  $\bar{\pi}_i:(G\cap A)\rightarrow
A_i/Y_i=\ov{A}_i$ be the projection corresponding to $\pi_i$.  Denote by
$\ov{K}=KY/Y$ the corresponding Carter subgroup of $\ov{G}=G/Y$.  Then
$\ov{G}$
satisfies conditions of the lemma. By
induction, $(\ov{K}\cap\ov{A})^{\bar{\pi}_i}$ is a Carter subgroup of
$\ov{A}_i$. Let $K_1$ be a complete preimage of $\ov{K}$ in $G$ and let $Q$ be
a Hall $p'$-subgroup of $K_1$. Then $(Q\cap A)^{\pi_i}$ is a Hall $p'$-subgroup
of $(K_1\cap A)^{\pi_i}$. In view of the proof of
\cite[Theorem~20.1.4]{KarMer}, we obtain that $K=N_{K_1}(Q)$ is a Carter
subgroup of $G$ and $(N_{K_1\cap A}(Q\cap A))^{\pi_i}$ is a Carter subgroup of
$A_i$. Thus we need to show that $(N_{K_1\cap A}(Q\cap A))^{\pi_i}=(N_{K_1\cap
S}(Q))^{\pi_i}$. By induction, the equality $(N_{\ov{K}\cap \ov{A}}(\ov{A}\cap
\ov{Q}))^{\bar{\pi}_i}=(N_{\ov{K}\cap \ov{G}}(\ov{Q}))^{\bar{\pi}_i}$ holds.
Thus we need to prove that $(N_Y(Q\cap A))^{\pi_i}=(N_Y(Q))^{\pi_i}$. Note also
that $(N_Y(Q\cap A))^{\pi_i}\leq N_{Y_i}((Q\cap A)^{\pi_i})$.

Since $S$ is a transitive and primitive nilpotent subgroup of $\Sym_k$, then
$k=r$ is prime and $S=\langle s\rangle$ is cyclic. If $r=p$, then $Q\cap A=Q$
and we have nothing to prove. Otherwise let $h$ be an $r$-element of $K$,
generating $S$ modulo $K\cap A$. Clearly $Q=(Q\cap A)\langle h\rangle$. Let
$t\in Y_i$
be an element of $N_{Y_i}((Q\cap A)^{\pi_i})$.
Then
$(t\cdot t^h\cdot\ldots\cdot t^{h^{r-1}})\in N_Y(Q)$ and
$t^{\pi_i}= (t\cdot t^h\cdot\ldots\cdot
t^{h^{r-1}})^{\pi_i}$, hence~${(N_Y(Q\cap A))^{\pi_i}\leq
N_{Y_i}((Q\cap A)^{\pi_i})\leq (N_Y(Q))^{\pi_i}\leq {(N_Y(Q\cap A))^{\pi_i}}}$.
\end{proof}

\begin{ttt}\label{CriterionOfExistence}
Let $G$ be a finite group.  Then $G$ contains a Carter subgroup if and only if
$G$ satisfies~{\em\bfseries (E)}.
\end{ttt}

\begin{proof}
We prove first the part ``only if''. 
Let $H$ be a minimal normal subgroup of $G$. Then  $H=T_1\times\ldots\times
T_k$, where $T_1\simeq\ldots\simeq T_k\simeq T$ are simple  groups. 

If $H$ is elementary Abelian (i.~e., $T$ is cyclic of prime order), then
$\Aut(T)$ is solvable and contains a Carter subgroup. Assume that $T$ is a
nonabelian simple group. Clearly $K$ is a Carter subgroup of $KH$. 
By \cite[Lemma~3]{Vdoconj} we obtain that $\Aut_{KH}(T_i)$ contains a Carter
subgroup for all $i$.

Now we prove the part ``if''. Again assume by contradiction that $G$ is a
counterexample of minimal order, i.~e., that $G$ does
not contain a Carter subgroup, but, $G$ satisfies {\bfseries (E)}. Let $H$ be a
minimal normal subgroup of $G$. Then $H=T_1\times\ldots\times
T_k$, where $T_1\simeq \ldots\simeq T_k\simeq T$, and $T$ is a 
finite simple group. 

By definition $G/H$ satisfies {\bfseries (E)}, thus, by induction, there exists
a Carter subgroup $\ov{K}$ of $\ov{G}=G/H$. Let $K$ be a complete preimage of
$\ov{K}$, then $K$ satisfies {\bfseries (E)}. If $K\not=G$, then, by induction
$K$ contains a Carter subgroup $K'$. Note that $K'$ is a Carter subgroup of
$G$. Indeed, assume that $x\in N_G(K')\setminus K'$. Since $K'H/H=\ov{K}$ is a
Carter subgroup of $\ov{G}$, we have that $x\in K$. But $K'$ is a Carter
subgroup of $K$, thus~${x\in K'}$. Hence $G=K$, i.~e. $G/H$ is nilpotent.

If $H$ is Abelian, then $G$ is solvable, therefore $G$ contains a Carter
subgroup. So assume that $T$ is a nonabelian finite simple group. We first
show that
$C_G(H)$ is trivial. Assume that $C_G(H)=M$ is nontrivial. Since $T$ is a
nonabelian simple group, it follows that $M\cap H=\{e\}$, so $M$ is nilpotent.
By Lemma
\ref{CarterSubgroupsInCompFactorsUnderHomomorphism} we have that $G/M$ satisfy 
{\bfseries (E)}. By induction we obtain that $G/M$
contains a Carter subgroup $\ov{K}$. Let $K'$ be a complete preimage of
$\ov{K}$ in $G$. Then $K'$ is solvable, hence contains a Carter subgroup $K$.
Like above we obtain that $K$ is a Carter subgroup of $G$, a
contradiction. Hence~${C_G(H)=\{e\}}$.

Since $H$ is a minimal normal subgroup of $G$, we obtain that
$\Aut_G(T_1)\simeq \Aut_G(T_2)\simeq\ldots\simeq \Aut_G(T_k)$.  Thus there
exists a monomorphism $$\varphi:G\rightarrow
\left(\Aut_G(T_1)\times\ldots\times\Aut_G(T_k)\right)\leftthreetimes
\Sym_k=G_1$$ and we identify $G$
with $G^\varphi$.  Denote by $K_i$ a Carter subgroup of $\Aut_G(T_i)$ and by
$A$ the subgroup $\Aut_G(T_1)\times\ldots\times\Aut_G(T_k)$. Since $G/H$ is 
nilpotent, then $K_iT_i=\Aut_G(T_i)$ and $G_1=\left(
K_1T_1\times\ldots\times K_kT_k\right)\leftthreetimes \Sym_k$. Let $\pi_i:G\cap
A\rightarrow (G\cap A)/C_{(G\cap A)}(T_i)$ be the canonical projections. Since
$G/(G\cap A)$ is transitive, we obtain that $(G\cap A)^{\pi_i}=K_iT_i$.

Since $\Aut_{G\cap A}(T_i)=K_iT_i$, hence $G\cap A$ satisfies
{\bfseries (E)}. By induction it contains a Carter subgroup $M$. By
\cite[Lemma~3]{Vdoconj} we obtain
that $M^{\pi_i}$ is a Carter subgroup of $K_iT_i$, therefore we may assume
$M^{\pi_i}=K_i$.  In particular, if  $R=(K_1\cap
T_1)\times\ldots\times(K_k\cap T_k)$, then  $M\leq N_G(R)$. 
 In view of
\cite{VdoAlmSimp}, Carter subgroups in every finite group are conjugate. Since
$(G\cap A)/H$ is nilpotent we obtain that $G\cap A=MH$, hence $G=N_G(M)H$.
Moreover $N_G(M)\cap A=M$, so $N_G(M)$ is solvable. Since $M$ normalizes $R$,
$M^{\pi_i}=K_i$, we obtain that $N_G(M)$ normalizes $R$,
hence $N_G(M)R$ is solvable. Therefore it contains a Carter subgroup $K$. By
Lemma~\ref{CarterInSubDirectProductOfSoluble}, $(K\cap A)^{\pi_i}$ is a Carter
subgroup of $(N_G(M)R\cap A)^{\pi_i}$ ($R$ plays the role of $T$ from
Lemma~\ref{CarterInSubDirectProductOfSoluble} in this case), so $(K\cap
A)^{\pi_i}=K_i$. Assume that $x\in N_G(K)\setminus K$. Since
$G/H=N_G(M)H/H=KH/H$ it follows that $x\in H$. Therefore 
$x^{\pi_i}\in (N_G(K)\cap A)^{\pi_i}\leq N_{T_i}((K\cap A)^{\pi_i})=K_i$. Since
$\bigcap_i \mathrm{Ker}(\pi_i)=\{e\}$, it follows that $x\in R\leq
N_G(M)R$. But $K$ is a Carter subgroup of $N_G(M)R$, hence $x\in K$.
This contradiction completes the proof.
\end{proof}

\section{Example}

In this section we construct an example showing, that we can not substitute
condition {\bfseries (E)} by a weaker condition: for every composition factor
$S$ of $G$, $\Aut_G(S)$ contains a Carter subgroup. This example also shows
that an extension of a group containing a Carter subgroup by a group containing
a Carter subgroup may fail to contain a Carter subgroup.

Consider $L=\Gamma SL_2(3^3)=PSL_2(3^3)\leftthreetimes \la\varphi\ra$,
where $\varphi$ is a field automorphism of $PSL_2(3^3)$. Let $X=(L_1\times
L_2)\leftthreetimes \Sym_2$, where $L_1\simeq L_2\simeq L$ and if
$\sigma=(1,2)\in\Sym_2\setminus\{e\}$, $(x,y)\in L_1\times L_2$, then
$\sigma(x,y)\sigma=(y,x)$ (permutation wreath product of $L$ and $\Sym_2$).
Denote
by $H=PSL_2(3^3)\times PSL_2(3^3)$ the minimal normal subgroup of $X$ and by
$M=L_1\times L_2$. Let
$G=(H\leftthreetimes \la(\varphi,\varphi^{-1})\ra)\leftthreetimes\Sym_2$ be a
subgroup of~$X$. Then the following statements hold:
\begin{itemize}
\item[1.] For every composition factor $S$ of $G$, $\Aut_G(S)$ contains a
Carter subgroup. 
\item[2.] $G\cap M\unlhd G$ contains a Carter subgroup.
\item[3.] $G/(G\cap M)$ is nilpotent.
\item[4.] $G$ does not contain a Carter subgroup.
\end{itemize}

1. Clearly we need to verify the statement for nonabelian composition factors
only. Every nonabelian composition factor $S$ of $G$ is isomorphic 
to $PSL_2(3^3)$ and $\Aut_G(S)=L$. In view of
\cite[Theorem~7.1]{VdoAlmSimp} we obtain that $L$ contains a Carter subgroup
(coinciding with a Sylow $3$-subgroup).

2. Since $(G\cap M)/H$ is nilpotent and from the previous statement we obtain
that $G\cap M$ satisfies {\bfseries (E)}, hence contains a Carter subgroup (it
is easy to see that a Sylow $3$-subgroup of $G\cap M$ is a Carter subgroup of
$G\cap M$).

3. Evident.

4. Assume that $K$ is a Carter subgroup of $G$. Then $KH/H$ is a Carter
subgroup of $G/H$. But $G/H$ is a nonabelian group of order $6$, hence
$G/H\simeq \Sym_3$ and $KH/H$ is a Sylow $2$-subgroup of $G/H$. In view
of \cite[Lemma~3]{Vdoconj} it follows that $\Aut_K(PSL_2(3^3))$ is a Carter
subgroup of $\Aut_{KH}(PSL_2(3^3))=PSL_2(3^3)$. But $PSL_2(3^3)$ does not
contain Carter subgroups in view of \cite[Theorem~7.1]{VdoAlmSimp}.
\vspace{1\baselineskip}

The author thanks Mazurov Vicktor Danilovoch for discussings on this paper,
that allow to improve the paper.

\end{document}